\documentclass[11pt,a4paper]{article}

 \usepackage{epsfig}
\usepackage{amssymb} 

\newcommand{\be}{\begin{equation}}
\newcommand{\ee}{\end{equation}}
\newcommand{\pd}{\partial}

\newenvironment{keywords}{ \noindent {\small\bf Key Words}:}{ }
\newenvironment{Acknowledgement}{ \noindent {\bf Acknowledgement}.}{ }

\begin{document}

\title{Software for Generation of Classes of Test Functions
with Known Local and Global Minima for Global Optimization}

\author{MARCO GAVIANO \\ Universit\`{a} di Cagliari, Italy
 \\ DMITRI E. KVASOV \\ Universit\`{a} di Roma ``La Sapienza'',
 Italy, and\\ University of Nizhni Novgorod, Russia
 \\ DANIELA LERA \\ Universit\`{a} di Cagliari, Italy\\
 YAROSLAV D. SERGEYEV\footnote{\scriptsize Corresponding author, e-mail:
yaro@si.deis.unical.it} \\ Universit\`{a} della Calabria, Italy,
and\\ University of Nizhni Novgorod, Russia}

 \maketitle

\begin{abstract}
A procedure for generating non-differentiable, continuously
differentiable, and twice continuously differentiable classes of
test functions for multiextremal multidimensional box-constrained
global optimization and a corresponding package of C subroutines
are presented. Each test class consists of 100 functions. Test
functions are generated by defining a convex quadratic function
systematically distorted by polynomials in order to introduce
local minima. To determine a class, the user defines the following
parameters: (i)~problem dimension, (ii)~number of local minima,
(iii)~value of the global minimum, (iv)~radius of the attraction
region of the global minimizer, (v)~distance from the global
minimizer to the vertex of the quadratic function. Then, all other
necessary parameters are generated randomly for all 100 functions
of the class. Full information about each test function including
locations and values of all local minima is supplied to the user.
Partial derivatives are also generated where possible.
\end{abstract}



\section{Introduction}
A wide literature is dedicated to development of numerical
algorithms for solving the global optimization problem (see, for
example, references given in \cite{Horst&Pardalos(1995)}). The
problem may be formulated as
 \be
  f^*=f(x^*)=\min_{x\in \mathcal{F} } f(x), \ \ \mathcal{F} \subset \mathbb{R}
  ^N,\label{glop}
 \ee
where $f(x)$ is a multiextremal and possibly non-differentiable
function and $\mathcal{F}$ is a compact set.

One of the approaches to studying and verifying validity of
numerical algorithms is their comparison on test problems (see,
e.g., \cite{Ali:et:al.(2003)}, \cite{Dixon&Szego(1978)},
\cite{Facchinei:et:al.(1997)}, \cite{Floudas&Pardalos(1990)},
\cite{Floudas:et:al(1999)}, \cite{Gaviano&Lera(1998)},
\cite{Horst&Pardalos(1995)}, \cite{Kalantari&Rosen(1986)},
\cite{Khoury:et:al.(1993)}, \cite{Li&Pardalos(1992)},
\cite{Locatelli(2003)}, \cite{More:et:al.(1981)},
\cite{Moshirvaziri(1994)}, \cite{Moshirvaziri:et:al(1996)},
 \cite{Pardalos(1987)}; \cite{Pardalos(1991)},
\cite{Pinter(2002)},  \cite{Schittkowski(1980)};
\cite{Schittkowski(1987)}, \cite{Schoen(1993)},
\cite{Sung&Rosen(1982)}). Many global optimization tests were
taken from real-life problems and for this reason comprehensive
information about them is not available. The number of local
minima may be unknown, as well as their locations, regions of
attraction, and even values (including that of the global
minimum).

Recently \cite{Gaviano&Lera(1998)} introduced two types of
functions with a priori known local minima and their regions of
attraction. The tests proposed take a convex quadratic function
(called hereafter `paraboloid') systematically distorted  by cubic
polynomials and by quintic polynomials to introduce local minima
and to construct test functions that are continuously
differentiable in some region $\Omega \supseteq \mathcal{F}$
(called hereafter `D-type' test functions) and twice continuously
differentiable in $\Omega \supseteq \mathcal{F}$ (called hereafter
`D2-type' test functions), where $\mathcal{F}$ is
from~(\ref{glop}) and $\Omega$ is a hyperrectangle.

To define a function of one of these types it is necessary to
determine a number of correlated parameters. Unfortunately, the
correlations do not allow simple and fast generation of the test
functions. Additionally, generation of different functions having
similar properties becomes difficult and non-intuitive when
dimension and/or number of local minima increase.

In this paper, in addition to the two types of test functions from
\cite{Gaviano&Lera(1998)}, the third type of non-differentiable
test functions (called hereafter `ND-type') is presented and a
generator for these three types of test functions is proposed. The
software to be introduced generates classes of test functions and
provides procedures for calculating the first order derivatives of
the D-type test functions and the first and second order
derivatives of the D2-type test functions.

Each class contains 100 functions and is defined by the following
parameters (the only ones to be determined by the user):
 \begin{enumerate}
  \item problem dimension;
  \item number of local minima;
  \item value of the global minimum;
  \item radius of the attraction region of the global minimizer;
  \item distance from the global minimizer to the vertex of the
  paraboloid.
 \end{enumerate}

The other necessary parameters (i.e., locations of all minimizers,
their regions of attraction, and values of minima) are chosen
randomly by the generator. After generation a special notebook
containing a complete description of all the functions from the
generated class is supplied to the user.

The rest of the paper is structured as follows. In
Section~\ref{sectionMatDescr}, a mathematical description of the
three types of test functions is given.
Section~\ref{sectionDescr}~introduces the generator and details of
its implementation. Section~\ref{sectionUsing}~is devoted to usage
of the generator.

\section{Mathematical description} \label{sectionMatDescr}

In this section, the three types of test functions are briefly
described. Let us start with the D-type and D2-type functions
(see~\cite{Gaviano&Lera(1998)}). A function $f(x)$ of the D-type
is determined over an admissible region $\Omega \supseteq
\mathcal{F}$, where $\mathcal{F}$ is from~(\ref{glop}) and
 \be
  \Omega=[a,b\,]=\{x \in \mathbb{R} ^N \, : \, a \leq x \leq
  b \}, \hspace{3mm} a < b, \hspace{2mm} a,b \in \mathbb{R}^N.
  \label{D}
 \ee
The function is constructed by modifying a paraboloid $Z$:
 \be
   Z:\ g(x)= \| x-T \| ^2 + \, t,\ \ x \in \Omega, \label{Z}
 \ee
(hereafter $\| \cdot \|$ denotes the Euclidean norm) with the
minimum $t$ at a point \mbox{$T \in {\rm int}(\Omega)$} in such a
way that the resulting function $f(x)$ has $m$, $m \geq 2$, local
minimizers: point~$T$ from (\ref{Z}) (we denote it by $M_1 := T$)
and points
 \be
  M_i \in {\rm int}(\Omega), \ \ M_i \neq T, \ M_i \neq M_j, \ \ i,j=2,\ldots,m,\ i \neq
  j. \label{M(i)}
 \ee
The paraboloid $Z$ from~(\ref{Z}) is modified by a function
$C_i(x)$, which is constructed by using cubic polynomials within
balls $S_i \subset \Omega$ around each point $M_i$,
$i=2,\ldots,m$, where
 \be
  S_i = \{x \in \mathbb{R} ^N : \, \| x - M_i \| \ \leq
  \rho_i,\, \rho_i > 0\}, \  \ i= 1, \ldots, m. \label{S(i)}
 \ee
Functions $Q_i(x)$, $i=2,\ldots,m$, use quintic polynomials to
determine the D2-type test functions.

Selection of radii $\rho_i$, $i=1,\ldots,m$, is carried out in
such a manner that sets $S_i$ from~(\ref{S(i)}) do not overlap:
 \be
  S_i \cap S_j = \varnothing,\ \ i,j=1,\ldots,m,\ i \neq j. \label{S(i)xS(j)}
 \ee
It is not required that each attraction region $S_i$, $i=1,
\ldots, m$, be entirely contained in $\Omega$. Note that we use
the notation ``attraction region'' with respect to the balls
$S_i$, $i=1,\ldots, m$, just for simplicity. Naturally, definition
of the real attraction region for each local minimizer will depend
on the method used for optimization and will change from one
algorithm to another.

Formally, D-type functions \cite{Gaviano&Lera(1998)} are described
as follows:
 \be
  f(x)=\left\{
   \begin{array}{ll}
     C_i(x), & x\in S_i, \, i \in \{2,\ldots,m\} , \\
     g(x), & x \notin S_2 \cup \ldots \cup S_m
     \, ,
   \end{array}
  \right. \label{f_cubic}
 \ee
where $g(x)$ is from~(\ref{Z}), sets $S_i$, $i=2,\ldots,m$,
from~(\ref{S(i)}) satisfy~(\ref{S(i)xS(j)}), and
 $$
  C_i(x) =\left ( \frac{2}{\rho_i^2} \frac{<\!x-M_i,\, T-M_i\!>}{\| x-M_i \| } - \frac{2}{\rho_i^3}A_i \right) \| x-M_i\| ^3 +
 $$
 \be
  + \left( 1- \frac{4}{\rho_i}\frac{<\!x-M_i,\, T-M_i\!>}{\|   x-M_i \| }+\frac{3}{\rho_i^2}A_i \right) \| x-M_i\| ^2 + f_i. \label{C_cubic}
 \ee
In~(\ref{C_cubic}) radii $\rho_i$, $i=2,\ldots,m$, determine the
sets $S_i$ from~(\ref{S(i)}), $<\! \cdot,\cdot \!>$ denotes the
usual scalar product, and the values $A_i$, $i=2,\ldots,m$, are
found as
 \be
  A_i= \| T-M_i\| ^2 + t - f_i, \label{A}
 \ee
where $f_1 = t$ and $f_i$, $i=2,\ldots,m$, are the function values
at local minimizers $M_i$:
 \be
  f_i = \min \{g(x) \, : x \in B_i\} - \gamma_i,\ \ \gamma_i > 0,
  \label{f(i)}
 \ee
where $B_i$ is the boundary of the ball $S_i$:
 \be
  B_i=\{x \in \mathbb{R} ^N : \, \| x - M_i \| \ = \rho_i, \, \rho_i > 0\}, \  \ i= 2, \ldots,
  m, \label{B(i)}
 \ee
and $\gamma_i$ is a parameter ensuring that the value $f_i$ is
less than the minimum of the paraboloid $Z$ from~(\ref{Z}) over
$B_i$.

Analogously, D2-type functions \cite{Gaviano&Lera(1998)} are
defined by
 \be
  f(x)=\left\{
   \begin{array}{ll}
     Q_i(x), & x\in S_i, \, i \in \{2,\ldots,m\} , \\
     g(x), & x \notin S_2 \cup \ldots \cup S_m
     \, ,
   \end{array}
  \right. \label{f_quintic}
 \ee
where
 \begin{center}
 $$
  Q_i(x) =\left [-\frac{6}{\rho_i^4} \frac{<\!x-M_i,\, T-M_i\!>}{\|   x-M_i \| } + \frac{6}{\rho_i^5}A_i +
  \frac{1}{\rho_i^3}(1-\frac{\delta}{2})\right]
  \| x-M_i\| ^5 +
 $$
 $$
  \left [\frac{16}{\rho_i^3} \frac{<\!x-M_i,\, T-M_i\!>}{\|   x-M_i \| } - \frac{15}{\rho_i^4}A_i -
  \frac{3}{\rho_i^2}(1-\frac{\delta}{2})\right]
  \| x-M_i\| ^4 +
 $$
 $$
  \left [-\frac{12}{\rho_i^2} \frac{<\!x-M_i,\, T-M_i\!>}{\|   x-M_i \| } + \frac{10}{\rho_i^3}A_i +
  \frac{3}{\rho_i}(1-\frac{\delta}{2})\right]
  \| x-M_i\| ^3 +
 $$
 \be
  \frac{1}{2}\delta \| x-M_i\| ^2 + f_i \label{Q_quintic}
 \ee
 \end{center}
with $A_i$ and $f_i$, $i=2,\ldots,m$, from~(\ref{A})
and~(\ref{f(i)}), and $\delta$ is an arbitrary positive real
number (see \cite[Lemma 3.1]{Gaviano&Lera(1998)}).

The properties of these functions have been studied by
\cite{Gaviano&Lera(1998)}. In particular, the following results
can be proved:

i. D-type functions~(\ref{f_cubic})--(\ref{C_cubic}) are
continuously differentiable in~$\Omega$~\cite[Lemma
2.1]{Gaviano&Lera(1998)}.

ii. D2-type functions~(\ref{f_quintic})--(\ref{Q_quintic}) are
twice continuously differentiable in~$\Omega$ \cite[Lemma
3.1]{Gaviano&Lera(1998)}.

Let us now describe the ND-type test functions, which are
continuous in $\Omega$ but non-differentiable in the whole region
$\Omega$. An analogous procedure is considered: the paraboloid $Z$
from~(\ref{Z}) is modified by a function $P_i(x)$ constructed from
second degree polynomials within each region $S_i \subset \Omega$
from~(\ref{S(i)}) in such a way that the resulting function $f(x)$
is continuous in the feasible region $\Omega$ from~(\ref{D}),
differentiable at each local minimizer $M_i$, $i=2,\ldots, m$,
from~(\ref{M(i)}), but generally non-differentiable at the points
of the boundaries $B_i$ of the balls $S_i$, $i=2, \ldots, m$,
determined by~(\ref{B(i)}). That is,
 \be
  f(x)=\left\{
   \begin{array}{ll}
     P_i(x), & x\in S_i, \, i \in \{2,\ldots,m\} , \\
     g(x), & x \notin S_2 \cup \ldots \cup S_m
     \, ,
   \end{array}
  \right. \label{f_nonsmooth}
 \ee
where $g(x)$ is from~(\ref{Z}), sets $S_i$, $i=2,\ldots,m$,
from~(\ref{S(i)}) satisfy~(\ref{S(i)xS(j)}), and
 \be
  P_i(x) =\left ( 1 -\frac{2}{\rho_i} \frac{<\!x-M_i,\, T-M_i\!>}{\|   x-M_i \| } + \frac{1}{\rho_i^2}A_i \right) \| x-M_i\| ^2 + f_i.
  \label{P_quadratic}
 \ee
In~(\ref{P_quadratic}) the values $\rho_i$, $A_i$, and $f_i$
($i=2,\ldots,m$) are determined in the same way as for the D- and
D2-type functions by formulae~(\ref{S(i)})--(\ref{S(i)xS(j)}),
(\ref{A}), and~(\ref{f(i)}), respectively.

\section{Generation of tests classes} \label{sectionDescr}

As one can see from the previous section, all three function types
have many parameters to be coordinated. Moreover, their
characteristics (for example, the mutual positions of the local
minimizers, the global minimizer, and the paraboloid vertex; the
size of the attraction regions of local minimizers; the function
values at local minima) influence the properties of the test
functions significantly from the point of view of global
optimization algorithms. For example, coincidence of the global
minimizer with the paraboloid vertex leads to generation of too
simple functions. Existence of many deep minima having narrow
regions of attraction can lead to the impossibility of global
minimizer location even by the most ``intelligent'' global
optimization algorithms. All these features should be added to the
general scheme from Section~\ref{sectionMatDescr} in order to
obtain well-structured test classes.

In the generator, the user sets just a few parameters defining a
desirable class while all the other parameters are chosen
randomly. The generator is also employed in maintaining conditions
distinguishing each class -- for example, the distance of the
global minimizer from the minimizer of the paraboloid, dependence
of the local minima values on the attraction regions sizes, etc.
Thus, the generator gives the researcher the ability to construct
classes of 100 test functions of arbitrary dimension with
arbitrary number of local minima.

This section describes how a class consisting of D-type test
functions is generated. Classes consisting of D2-type and ND-type
functions are constructed analogously.

Each test class generated by the introduced software contains 100
test functions $f(x)$ and is defined by the following parameters
to be fixed by the user:
\begin{enumerate}
 \item the problem dimension $N$, $N \geq 2$;
 \item the number of local minimizers $m$, $m \geq 2$, including the
 minimizer $T$ for the paraboloid~(\ref{Z}) (all the minimizers
 are chosen randomly);
 \item the global minimum value $f^*$, the same for all the
 functions of the class;
 \item the radius $\rho ^*$ of the attraction region of the global
 minimizer $x^*$;
 \item the distance $r^*$ from the paraboloid vertex $T$ to the global
 minimizer \mbox{$x^*\in \Omega$} (whose coordinates are also chosen
 randomly).
\end{enumerate}
By changing these parameters the user can create classes with
different properties.

Each function of a test class is specified by its number~$n$, $1
\leq n \leq 100$. The other parameters of the functions
from~(\ref{Z})--(\ref{P_quadratic}) are chosen randomly by means
of the random number generator proposed in~\cite{Knuth(1997)}.

The input parameters $f^*$, $r^*$, and $\rho ^*$ must be chosen
in such a way that the following simple conditions are satisfied:
 \be
  f^* < t \label{f*}
 \ee
(which means that the global minimizer is not a vertex of the
paraboloid; this requirement allows us to avoid too simple
functions with a global minimum at the vertex of the paraboloid
$Z$ from~(\ref{Z})),
 \be
  0 < r ^* < 0.5\min_{1 \leq j \leq N} |  b(j) - a(j) |  \label{r_global_cond}
 \ee
(i.e., the global minimizer $x^*$ belongs to the admissible region
$\Omega$ even in the case when the paraboloid vertex $T$ is at the
center of $\Omega$), and
 \be
  0 < \rho ^* \leq 0.5r^*. \label{rho_global}
 \ee

Note that it is not required that each attraction region $S_i$,
$i=1, \ldots, m$, from~(\ref{S(i)}) entirely belongs to $\Omega$.

The admissible region $\Omega$ is taken as $\Omega=[-1,1]^N$ and
the minimal value of the paraboloid~(\ref{Z}) is fixed at $t=0$ by
default (naturally, these parameters can be changed by the user).

Let us discuss in more detail the random procedure generating
parameters for test functions. (The unique difference for the
D2-type is that the parameter $\delta$ from~(\ref{Q_quintic}) is
required; this parameter is chosen randomly from the open interval
$(0,\Delta)$, where $\Delta$ is a positive number taken by default
$\Delta=10$.) Hereafter the vertex $T$ from~(\ref{Z}) in the
set~(\ref{M(i)}) of local minimizers has the index 1, $M_1 := T$,
and the global minimizer $x^*$ has the index 2, $M_2 := x^*$.
Naturally, among the minimizers $M_i$, $i=3, \ldots, m$, another
global minimizer $y^* \neq x^*$ can be generated.

First, coordinates of the paraboloid vertex $T$, coordinates of
the global minimizer~$x^*$, and coordinates of the remaining local
minimizers $M_i$ (controlling the satisfaction of~(\ref{M(i)}))
are chosen randomly. Then, the attraction regions radii $\rho_i$,
$i \neq 2$, from~(\ref{S(i)}) are determined: to do this the
attraction regions of each local minimizer from~(\ref{M(i)}) ($i
\neq 2$ because the attraction region of the global minimizer is
fixed: $\rho_2 = \rho^*$) are expanded until
condition~(\ref{S(i)xS(j)}) is not violated. Finally, values of
the function $f(x)$ at local minima $M_i$, $i=3, \ldots, m$, are
fixed by choosing random values $\gamma_i$, $i=3, \ldots, m$,
from~(\ref{f(i)}) (recall that $f_1=t$ and $f_2 = f^*$).

Let us consider these three principal operations in detail.

Coordinates of the local minimizers $M_i$, $i=3, \ldots, m$,
from~(\ref{M(i)}), coordinates of the vertex $T$ of the
paraboloid~(\ref{Z}), and location of the global minimizer $x^*$
are chosen randomly at the intersection of $\Omega$ and the sphere
of radius $r^*$ with a center at~$T$ so that~(\ref{M(i)}) is
satisfied. For the positioning of $x^*$ we use generalized
spherical coordinates
 $$
  x^*_j := T_j + r^* \cos \phi _j \prod_{k=1}^{j-1}\sin \phi _k, \ \
  j=1,\ldots,N-1,
 $$
 \be
  x^*_N := T_N + r^* \prod_{k=1}^{N-1}\sin \phi _k, \label{x_polar}
 \ee
where the components of the vector
 $$
  \phi = (\phi _1, \ldots, \phi _N) \in \Phi = \{0 \leq \phi_1 \leq \pi;
  \ 0 \leq \phi_j \leq 2\pi, \ j=2,\ldots,N \}
 $$
are chosen randomly. In this case, if some $x^*_k \notin \Omega$,
$1 \leq k \leq N$, this coordinate is redefined as
 $$
  x^*_k := 2 T_k - x^*_k.
 $$

After selection of coordinates of the paraboloid vertex $T$ and of
the global minimizer $x^*$, coordinates of the points $M_i$, $i=3,
\ldots, m$, are generated in such a way that beside
condition~(\ref{M(i)}) the condition
 \be
  \| M_i - x^* \| - \rho ^* = \zeta, \ \ \zeta > 0
  \label{gap}
 \ee
is satisfied with some positive parameter $\zeta$. This condition
follows from~(\ref{S(i)xS(j)}) and does not allow the local
minimizers to be very close to the attraction region of the global
minimizer $x^*$. Thus, in~(\ref{gap}) the parameter $\zeta$ should
not be too small. The value $\zeta = \rho^*$ is chosen by default.

The next step of the test function construction sets attraction
regions. Each value $\rho_i$, $i \neq 2$, from~(\ref{S(i)}) is
initially calculated as half of the minimum distance between the
minimizer $M_i$ and the remaining local minimizers
 $$
  \rho_i := 0.5 \min_{1 \leq j \leq m,\ j \neq i} \|   M_i-M_j\| , \ \ i=1,\ldots,m, \ i \neq 2,
 $$
 \be \label{rho(i)init}
  \rho_2 := \rho ^*
 \ee
(in such a way that the attraction regions from~(\ref{S(i)}) do
not overlap). Then, an attempt to increase the values $\rho_i$,
$i=1,\ldots,m$, $i \neq 2$ (i.e., an attempt to enlarge the
attraction regions) is made:
 \be
  \rho _i := \max \left( \rho_i,\min_{1 \leq j \leq m, \ j \neq
  i} \{ \| M_i-M_j\| - \rho_j \}\right), \ \
  i=1,\ldots,m, \ i \neq 2. \label{rho(i)expand}
 \ee
Because of the recursive character of
formulae~(\ref{rho(i)expand}), an expansion of the attraction
regions depends on the order in which these regions are selected
(an ascending order of the indices is chosen).

Finally, the values of the radii $\rho_i$ are corrected by the
weight coefficients $w_i$:
 $$
  \rho_i := w_i\,\rho_i, \ \ i=1, \ldots, m,
 $$
where $0 < w_i \leq 1$, $i=1,\ldots, m$, and the values $w_i$ are
chosen by default as
 \be
  w_i = 0.99, \ \ i=1, \ldots, m, \ i \neq 2, \hspace{7mm} {\rm
  and} \hspace{7mm} w_2 = 1.  \label{rho(i)weight}
 \ee

At the last step the function values $f_i$, $i=3, \ldots, m$, at
the local minima are generated by using formula~(\ref{f(i)}),
where $\gamma_i$ must be specified. Each value $\gamma_i$, $i=3,
\ldots, m$, is chosen (note that the values $\gamma_1$ and
$\gamma_2$ are not considered because the function values $f_1=t$
at the paraboloid vertex and $f_2=f^*$ at the global minimizer
have been fixed by the user without using~(\ref{f(i)})) as the
minimum of two values generated randomly from the open intervals
$(\rho_i, 2\rho_i)$ and $(0, Z_{B_i} - f^*)$, where $Z_{B_i}$ is
the minimum of the paraboloid $Z$ from~(\ref{Z}) over $B_i$
from~(\ref{B(i)}). In such a way, the values $f_i$ in~(\ref{f(i)})
depend on radii $\rho_i$ of the attraction regions $S_i$, $i=3,
\ldots, m$, and at the same time the following condition is
satisfied:
 $$
  f^* \leq f_i,\ \ i=3, \ldots, m.
 $$
Note that dependence of the function values at local minima on the
radii of the attraction regions is not respected by the global
optimum value $f_2=f^*$ because the user defines the function
value at the global minimizer and the radius $\rho ^*$ of its
region of attraction directly when choosing the corresponding test
class.

\begin{figure}[t]
\centerline{\psfig{file=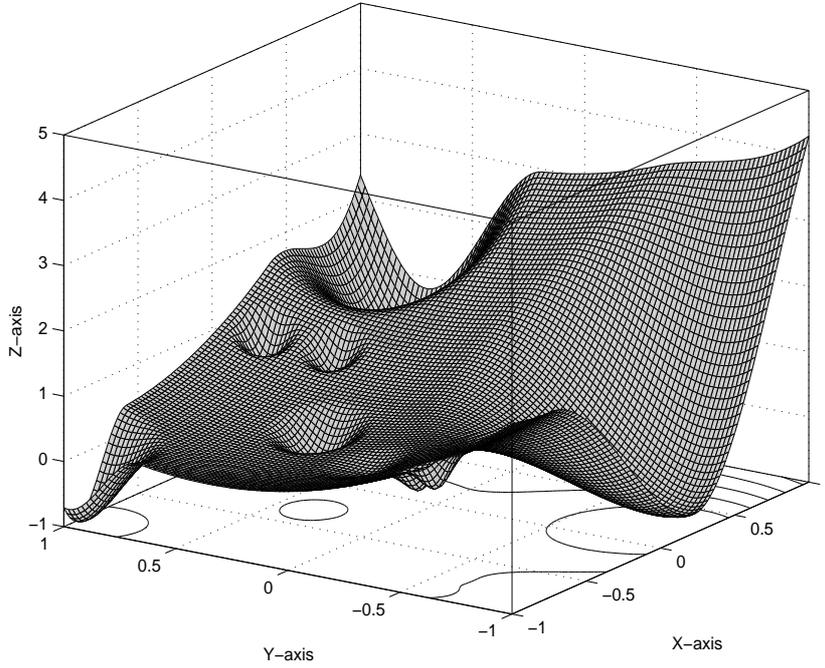,width=110mm,height=90mm,angle=0,silent=}}
\caption{The function number~9 from a class of two-dimensional
D-type test functions with~10 local minima} \label{fig1}
\end{figure}

Figure~\ref{fig1} shows an example of the D-type test function.
This function is defined in the region $\Omega=[-1,1]^2$ and is
number~9 in the class of D-type functions with the following
parameters:
 \begin{enumerate}
  \item dimension $N=2$;
  \item number of local minima $m=10$;
  \item value of the global minimum $f^*=-1$;
  \item radius of the attraction region of the global minimizer $\rho^*=\frac{1}{3}$;
  \item distance from the global minimizer $x^*$ to the vertex $T$ of the
  paraboloid from~(\ref{Z}) is $r^*=\frac{2}{3}$.
 \end{enumerate}
The generated global minimizer of this function is $x^*=(-0.911,
0.989)$ and the paraboloid minimizer is $T=(-0.711, 0.353)$.

\section{Usage of the test classes generator} \label{sectionUsing}

The generator package has been written in ANSI Standard C and
successfully tested on Windows and UNIX platforms. Our
implementation follows the procedure described in
Section~\ref{sectionDescr}. First, the general structure of the
package is described, then instructions for using the test classes
generator (called hereafter GKLS-generator) are given.

\subsection{Structure of the package}

The package includes the following files:
 \begin{description}
  \item[\bf gkls.c] -- the main file;
  \item[\bf gkls.h] -- the header file that users should include in
  their application projects in order to call subroutines from
  the file {\bf gkls.c};
  \item[\bf rnd\_gen.c] -- the file containing the uniform random number
  generator proposed in Knuth~[\cite{Knuth(1997)};
  \cite{Knuth:HomePage}];
  \item[\bf rnd\_gen.h] -- the header file for linkage to the file {\bf
  rnd\_gen.c};
  \item[\bf example.c] -- an example of the GKLS-generator usage;
  \item[\bf Makefile] -- an example of a UNIX makefile provided to
  UNIX users for a simple compilation and linkage of separate files
  of the application project.
 \end{description}

For implementation details the user can consult the C codes. Note
that the random number generator in {\bf rnd\_gen.c} uses the
logical-and operation `\&' for efficiency, so it is not strictly
portable unless the computer uses two's complement representation
for integer. It does not limit portability of the package because
almost all modern computers are based on two's complement
arithmetic.

\subsection{Calling sequence for generation and usage of the tests
classes}

Here we describe how to generate and use classes of the ND-, D-,
and D2-type test functions. Again, we concentrate on the D-type
functions. The operations for the remaining two types are
analogous.

To utilize the GKLS-generator the user must perform the following
steps:
\begin{description}
 \item[\bf Step 1.] Input of the parameters defining a specific test class.
 \item[\bf Step 2.] Generating a specific test function of the defined test
 class.
 \item[\bf Step 3.] Evaluation of the generated test function and, if
 necessary, its partial derivatives.
 \item[\bf Step 4.] Memory deallocating.
\end{description}
Let us consider these steps in turn.

\subsubsection{Input of the parameters defining a specific test
class} \label{sectionInput}

This step is subdivided into: (a) defining the parameters of the
test class, (b) defining the admissible region $\Omega$, and (c)
checking (if necessary).

-- {\it (a) Defining the parameters of the test class}. The
parameters to be defined by the user determine a specific class
(of the ND-, D- or D2-type) of 100 test functions (a specific
function is retrieved by its number). There are the following
parameters:
 \begin{description}
  \item[{\it GKLS\_dim}] -- ({\bf unsigned int}) dimension $N$
  (from~(\ref{glop})) of test functions; $N \geq 2$ (since
  multidimensional problems are considered in~(\ref{glop})) and
  $N<$ NUM\_RND in {\bf rnd\_gen.h}; this value is limited by the
  power of \mbox{{\bf unsigned int}}-representation; default
  $N=2$;
  \item[{\it GKLS\_num\_minima}] -- ({\bf unsigned int}) number $m$
  (from~(\ref{M(i)})) of local minima including the paraboloid $Z$
  minimum (from~(\ref{Z})) and the global minimum; $m \geq 2$; the
  upper bound of this parameter is limited by the power of
  \mbox{{\bf unsigned int}}-representation;
  default $m=10$;
  \item[{\it GKLS\_global\_value}] -- ({\bf double}) global minimum value $f^*$ of
  $f(x)$; condition~(\ref{f*}) must be satisfied; the default
  value is $-1.0$ (defined in the file {\bf gkls.h} as a constant~GKLS\_GLOBAL\_MIN\_VALUE);

  \item[{\it GKLS\_global\_dist}] -- ({\bf double}) distance $r^*$ from the paraboloid
  vertex $T$ in~(\ref{Z}) to the global minimizer $x^* \in \Omega$
  of $f(x)$; condition~(\ref{r_global_cond}) must be
  satisfied; the default value is
  $$
    {\rm GKLS\_global\_dist}\ \stackrel{\rm def}{=} \ \min_{1 \leq j
    \leq N} |  b(j) - a(j) |  /\, 3,
  $$
  where the vectors $a$ and $b$ determine the admissible region
  $\Omega$ in~(\ref{D});
  \item[{\it GKLS\_global\_radius}] -- ({\bf double}) radius $\rho ^*$ of the attraction
  region of the global
  minimizer $x^*\in \Omega$ of $f(x)$; condition~(\ref{rho_global}) must be
  satisfied; the default value is
  $$
    {\rm GKLS\_global\_radius}\ \stackrel{\rm def}{=} \ \min_{1 \leq
    j \leq N} |  b(j) - a(j) |  /\, 6.
  $$
 \end{description}
The user may call subroutine \mbox{{\it GKLS\_set\_default}()} to
set the default values of these five variables.

-- {\it (b) Defining the admissible region} $\Omega$. With $N$
determined, the user must allocate dynamic arrays {\it
GKLS\_domain\_left} and {\it GKLS\_domain\_right} to define the
boundary of the hyperrectangle $\Omega$. This is done by calling
subroutine

{\bf int} {\it GKLS\_domain\_alloc} ();\\
which has no parameters and returns the following error codes
defined in~{\bf gkls.h}:
\begin{description}
 \item[\bf GKLS\_OK] -- no errors;
 \item[\bf GKLS\_DIM\_ERROR] -- the problem dimension is out of range;
 it must be greater than or equal to 2 and less than NUM\_RND defined
 in {\bf rnd\_gen.h};
 \item[\bf GKLS\_MEMORY\_ERROR] -- there is not enough memory to
 allocate.
\end{description}

The same subroutine defines the admissible region $\Omega$. The
default value $\Omega =[-1,1]^N$ is set by \mbox{{\it
GKLS\_set\_default}()}.

-- {\it (c) Checking}. The following subroutine allows the user to
check validity of the input parameters:

{\bf int} {\it GKLS\_parameters\_check} ().\\
It has no parameters and returns the following error codes (see
{\bf gkls.h}):
 \begin{description}
  \item[\bf GKLS\_OK] -- no errors;
  \item[\bf GKLS\_DIM\_ERROR] -- problem dimension error;
  \item[\bf GKLS\_NUM\_MINIMA\_ERROR] -- number of local minima error;
  \item[\bf GKLS\_BOUNDARY\_ERROR] -- the admissible region boundary
  vectors are ill-defined;
  \item[\bf GKLS\_GLOBAL\_MIN\_VALUE\_ERROR] -- the global minimum
  value is not less than the paraboloid~(\ref{Z}) minimum value $t$
  defined in {\bf gkls.h} as a constant GKLS\_PARABOLOID\_MIN;
  \item[\bf GKLS\_GLOBAL\_DIST\_ERROR] -- the parameter $r^*$ does not
  satisfy~(\ref{r_global_cond});
  \item[\bf GKLS\_GLOBAL\_RADIUS\_ERROR] -- the parameter $\rho^*$ does not
  satisfy~(\ref{rho_global}).
 \end{description}

\subsubsection{Generating a specific test function of the defined
test class}

After a specific test class has been chosen (i.e., the input
parameters have been determined) the user can generate a specific
function that belongs to the chosen class of 100 test functions.
This is done by calling subroutine

{\bf int} {\it GKLS\_arg\_generate} ({\bf unsigned int} {\it nf}); \\
where
\begin{description}
 \item[\it nf] -- the number
 of a function from the test class (from~1~to~100).
\end{description}
This subroutine initializes the random number generator, checks
the input parameters, allocates dynamic arrays, and generates a
test function following the procedure of
Section~\ref{sectionDescr}. It returns an error code that can be
the same as for subroutines {\it GKLS\_parameters\_check}() and
{\it GKLS\_domain\_alloc}(), or additionally:
\begin{description}
 \item[\bf GKLS\_FUNC\_NUMBER\_ERROR] -- the number of a test function to
 generate exceeds 100 or it is less than 1.
\end{description}

{\it GKLS\_arg\_generate}() generates the list of all local minima
and the list of the global minima as parts of the structures {\it
GKLS\_minima} and {\it GKLS\_glob}, respectively. The first
structure gathers the following information about all local minima
(including the paraboloid minimum and the global one): coordinates
of local minimizers, local minima values, and attraction regions
radii. The second structure contains information about the number
of global minimizers and their indices in the set of local
minimizers. It has the following fields:
 \begin{description}
  \item[\it num\_global\_minima] -- ({\bf unsigned int}) total number of
  global minima;
  \item[\it gm\_index] -- ({\bf unsigned int *}) list of indices of generated
  minimizers, which are the global ones (elements 0 to ($\mbox{\it
  num\_global\_minima} - 1$) of the list) and the local ones
  (the remaining elements of the list).
 \end{description}
The elements of the list {\it GKLS\_glob}.{\it gm\_index} are
indices to a specific minimizer in the first structure {\it
GKLS\_minima} characterized by the following fields:
 \begin{description}
  \item[\it local\_min] -- ({\bf double **}) list of local minimizers
  coordinates;
  \item[\it f] -- ({\bf double *}) list of local minima values;
  \item[\it rho] -- ({\bf double *}) list of attraction regions radii;
  \item[\it peak] -- ({\bf double *}) list of parameters $\gamma_i$ values
  from~(\ref{f(i)});
  \item[\it w\_rho] -- ({\bf double *}) list of parameters $w_i$ values
  from~(\ref{rho(i)weight}).
 \end{description}
The fields of these structures can be useful if one needs to study
properties of a specific generated test function more deeply.

\subsubsection {Evaluation of a generated test function or its partial
derivatives}

While there exists a structure {\it GKLS\_minima} of local minima,
the user can evaluate a test function (or partial derivatives of
D- and D2-type functions) that is determined by its number (a
parameter to the subroutine {\it GKLS\_arg\_generate}()) within
the chosen test class. If the user wishes to evaluate another
function within the same class he should deallocate dynamic arrays
(see the next subsection) and recall the generator {\it
GKLS\_arg\_generate}() (passing it the corresponding function
number) without resetting the input class parameters (see
subsection~\ref{sectionInput}). If the user wishes to change the
test class properties he should reset also the input class
parameters.

Evaluation of an ND-type function is done by calling subroutine

{\bf double} {\it GKLS\_ND\_func} ({\it x}).\\
Evaluation of a D-type function is done by calling subroutine

{\bf double} {\it GKLS\_D\_func} ({\it x}).\\
Evaluation of a D2-type function is done by calling subroutine

{\bf double} {\it GKLS\_D2\_func} ({\it x}).\\
All these subroutines have only one input parameter
\begin{description}
 \item[\it x] -- ({\bf double *}) a point $x \in \mathbb{R} ^N$
 where the function must be evaluated.
\end{description}
All the subroutines return a test function value corresponding to
the point $x$. They return the value GKLS\_MAX\_VALUE (defined in
{\bf gkls.h}) in two cases: (a) vector~$x$ does not belong to the
admissible region $\Omega$ and (b) the user tries to call the
subroutines without generating a test function.

The following subroutines are provided for calculating the partial
derivatives of the test functions (see Appendix).

Evaluation of the first order partial derivative of the D-type
test functions with respect to the variable $x_j$
(see~(\ref{df_cubic})--(\ref{dC_cubic}) in Appendix) is done by
calling subroutine

{\bf double} {\it GKLS\_D\_deriv} ({\it j}, {\it x}). \\
Evaluation of the first order partial derivative of the D2-type
test functions with respect to the variable $x_j$
(see~(\ref{df_quintic})--(\ref{dQ_quintic}) in Appendix) is done
by calling subroutine

{\bf double} {\it GKLS\_D2\_deriv1} ({\it j}, {\it x}). \\
Evaluation of the second order partial derivative of the D2-type
test functions with respect to the variables $x_j$ and $x_k$ (see
in Appendix the
formulae~(\ref{d2(jk)f_quintic})--(\ref{d2(jk)Q_quintic}) for the
case $j \neq k$
and~(\ref{d2(jj)f_quintic})--(\ref{d2(jj)Q_quintic}) for the case
$j = k$) is done by calling subroutine

{\bf double} {\it GKLS\_D2\_deriv2} ({\it j}, {\it k}, {\it x}). \\
Input parameters for these three subroutines are:
\begin{description}
 \item[\it j, k] -- ({\bf unsigned int}) indices of the variables
 (that must be in the range from 1 to {\it GKLS\_dim}) with respect to
 which the partial derivative is evaluated;
 \item[\it x] -- ({\bf double *}) a point $x \in \mathbb{R} ^N$ where
 the derivative must be evaluated.
\end{description}
All subroutines return the value of a specific partial derivative
corresponding to the point $x$ and to the given direction. They
return the value GKLS\_MAX\_VALUE (defined in {\bf gkls.h}) in
three cases: (a) index ($j$ or $k$) of a variable is out of the
range [1,{\it GKLS\_dim}]; (b) vector $x$ does not belong to the
admissible region~$\Omega$; (c) the user tries to call the
subroutines without generating a test function.

Subroutines for calculating the gradients of the D- and D2-type
test functions and for calculating the Hessian matrix of the
D2-type test functions at a given feasible point are also
provided. These are

{\bf int} {\it GKLS\_D\_gradient} ({\it x}, {\it g}), \\

{\bf int} {\it GKLS\_D2\_gradient} ({\it x}, {\it g}), \\

{\bf int} {\it GKLS\_D2\_hessian} ({\it x}, {\it h}). \\
Here
\begin{description}
 \item[\it x] -- ({\bf double *}) a point $x \in \mathbb{R} ^N$
 where the gradient or Hessian matrix must be evaluated;
 \item[\it g] -- ({\bf double *}) a pointer to the gradient vector
 calculated at {\it x};
 \item[\it h] -- ({\bf double **}) a pointer to the Hessian matrix
 calculated at {\it x}.
\end{description}
Note that before calling these subroutines the user must allocate
dynamic memory for the gradient vector {\it g} or the Hessian
matrix {\it h} and pass the pointers {\it g} or {\it h} as
parameters of the subroutines.

These subroutines call the subroutines described above for
calculating the partial derivatives and return an error code ({\bf
GKLS\_DERIV\_EVAL\_ERROR} in the case of an error during
evaluation of a particular component of the gradient or the
Hessian matrix, or {\bf GKLS\_OK} if there are no errors).

\subsubsection {Memory deallocating}

When the user concludes his work with a test function he should
deallocate dynamic arrays allocated by the generator. This is done
by calling subroutine

{\bf void} {\it GKLS\_free} ({\bf void});\\
with no parameters.

When the user abandons the test class he should deallocate dynamic
boundaries vectors {\it GKLS\_domain\_left} and {\it
GKLS\_domain\_right} by calling subroutine

{\bf void} {\it GKLS\_domain\_free} ({\bf void});\\
again with no parameters.

It should be finally highlighted that if the user, after
deallocating memory, wishes to return to the same class,
generation of the class with the same parameters produces the same
100 test functions.

An example of the generation and use of some of the test classes
can be found in the file {\bf example.c}.

\vspace{5mm}
\begin{Acknowledgement}
This research was partially supported by the following projects:
FIRB RBAU01JYPN, FIRB RBNE01WBBB, and RFBR 01-01-00587. The
authors thank Associate Editor Michael Saunders and an anonymous
referee for their subtle suggestions.
\end{Acknowledgement}

\appendix
\section{APPENDIX}  
\subsection*{Formulae of derivatives of the D- and D2-type test
functions}

\renewcommand{\theequation}{A.\arabic{equation}}
\setcounter{equation}{0}  

In this section, analytical expressions of the partial derivatives
of the D- and D2-type test functions are given. We denote by
$T=(T_1, \ldots, T_N)$ the minimizer of the paraboloid $Z$
from~(\ref{Z}) and by $M_i=(m^i_1,...,m^i_N)$, $i=2,...m$, the
local minima (from~(\ref{M(i)})) of a test function. Thus, for a
D-type test function $f(x)$ given
by~(\ref{f_cubic})--(\ref{C_cubic}) we have
(see~\cite{Gaviano&Lera(1998)}):
 \be
  \frac{\pd f(x)}{\pd x_j}=\left\{
   \begin{array}{ll}
     \frac{\pd C_i (x)}{\pd x_j}, & x\in S_i, \, i \in \{2,\ldots,m\} , \\
     2(x_j-T_j), & x \notin S_2 \cup \ldots \cup S_m
     \, ,
   \end{array}
  \right. \label{df_cubic}
 \ee
for $j = 1,\ldots,N$, and
\begin{eqnarray} \label{dC_cubic}
  \frac{\partial C_i (x)}{\partial x_j} & =  & \frac2{\rho_i^2}
  h_j(x) \|x-M_i\|+3 \left( \frac 2 {\rho_i^2}
  \frac{<\!x-M_i,T-M_i\!>}{\| x-M_i \|}-
  \frac2{\rho_i^3}A_i \right) \times \nonumber \\
  &  & \hspace{-5mm} \times \; (x_j-m^i_j) \| x-M_i \| - \frac4{\rho_i} h_j(x) + \\
  &  & \hspace{-5mm} + \; 2\left(1- \frac4{\rho_i} \frac{<\!x-M_i,T-M_i\!>} {\| x-M_i \|}
  + \frac3{\rho_i^2}A_i \right) (x_j-m^i_j), \nonumber
\end{eqnarray}
\vspace{6pt} \noindent with $ h_j(x)=(T_j-m^i_j) \|x-M_i\|
-<\!x-M_i,T-M_i\!> (x_j-m^i_j)/\|x-M_i\| $.

The first order partial derivatives of the D2-type test functions
$f(x)$ given by \mbox{(\ref{f_quintic})--(\ref{Q_quintic})} are
calculated as follows (see~\cite{Gaviano&Lera(1998)}):
 \be
  \frac{\pd f(x)}{\pd x_j}=\left\{
   \begin{array}{ll}
     \frac{\pd Q_i (x)}{\pd x_j}, & x\in S_i, \, i \in \{2,\ldots,m\} , \\
     2(x_j-T_j), & x \notin S_2 \cup \ldots \cup S_m
     \, ,
   \end{array}
  \right. \label{df_quintic}
 \ee

for $j = 1,\ldots, N$, and
\begin{eqnarray} \label{dQ_quintic}
 \frac{\pd Q_i (x)}{\pd x_j} & = & -\frac6{\rho_i^4} h_j(x)
 \|x-M_i\|^3 + 5 (x_j-m^i_j) \| x-M_i \|^3 \times \nonumber \\
 &  & \hspace{-5mm} \times \left[ -\frac{6}{\rho_i^4} \frac{<\!x-M_i,T-M_i\!>}{\| x-M_i
 \|}+\frac{6}{\rho_i^5}A_i + \frac1{\rho_i^3} (1-\frac{\delta}2 )
 \right] + \nonumber \\
 &  & \hspace{-5mm} + \; \frac{16}{\rho_i^3} h_j(x) \|x-M_i\|^2 + 4 (x_j-m^i_j) \|
 x-M_i \|^2 \times \nonumber \\
 &  & \hspace{-5mm} \times \left[ \frac{16}{\rho_i^3} \frac{<\!x-M_i,T-M_i\!>}{\| x-M_i
 \|}-\frac{15}{\rho_i^4}A_i - \frac3{\rho_i^2} (1-\frac{\delta}2 )
 \right] - \\
 &  & \hspace{-5mm} -\; \frac{12}{\rho_i^2} h_j(x) \|x-M_i\| + 3 (x_j-m^i_j) \| x-M_i \|
 \times \nonumber \\
 &  & \hspace{-5mm} \times \left[ -\frac{12}{\rho_i^2} \frac{<\!x-M_i,T-M_i\!>}{\| x-M_i
 \|}+\frac{10}{\rho_i^3} A_i + \frac3{\rho_i} (1-\frac{\delta}2 )
 \right] + \nonumber \\
 &  & \hspace{-5mm} + \; \delta (x_j-m^i_j), \nonumber
\end{eqnarray}
\vspace{6pt} \noindent with $ h_j(x)=(T_j-m^i_j) \|x-M_i\|
-<\!x-M_i,T-M_i\!> (x_j-m^i_j)/\|x-M_i\|$.

Let us now consider the second order derivatives $\pd ^2f(x) /\pd
x_j\pd x_k$ and $\pd ^2 f(x)/\pd x_j^2$ of the D2-type test
functions $f(x)$. For mixed partial derivatives $\pd ^2f(x) /\pd
x_j\pd x_k$ we have
 \be \label{d2(jk)f_quintic}
  \frac{\pd ^2f(x)}{\pd x_j\pd x_k}= \left\{
  \begin{array}{ll}
     \frac{\pd ^2Q_i (x)}{\pd x_j\pd x_k}, & x\in S_i, \, i \in \{2,\ldots,m\} , \\
     0, & x \notin S_2 \cup \ldots \cup S_m \, ,
   \end{array} \right.
 \ee
for $j,k=1, \ldots, N, \ j\neq k$, and
 \begin{eqnarray} \label{d2(jk)Q_quintic}
  \frac{\pd^2 Q_i (x)}{\pd x_j\pd x_k} & = & -\frac6{\rho_i^4}
  \left[\frac{\pd h_j(x)}{\pd x_k} \|x-M_i\|^3+3h_j(x)(x_k-m^i_k)
  \|x-M_i\| \right] - \nonumber \\
  &  & \hspace{-5mm} -\; \frac{30}{\rho_i^4}h_k(x) (x_j-m^i_j) \|x-M_i\| +
  15 (x_j-m^i_j)(x_k-m^i_k)\|x-M_i\| \times \nonumber \\
  &  & \hspace{-5mm} \times \left[ -\frac6{\rho_i^4} \frac{<\!x-M_i,T-M_i\!>}{\| x-M_i \|}+
  \frac{6}{\rho_i^5}A_i +\frac1{\rho_i^3} (1-\frac{\delta}2 ) \right] + \nonumber \\
  &  & \hspace{-5mm} + \; \frac{16}{\rho_i^3} \left[\frac{\pd h_j(x)}{\pd x_k}
  \|x-M_i\|^2+2h_j(x)(x_k-m^i_k) \right] + \nonumber \\
  &  & \hspace{-5mm} + \; \frac{64}{\rho_i^3}h_k(x)(x_j-m^i_j) + 8 (x_j-m^i_j)(x_k-m^i_k) \times \nonumber \\
  &  & \hspace{-5mm} \times \left[ \frac{16}{\rho_i^3} \frac{<\!x-M_i,T-M_i\!>}{\| x-M_i \|}-
  \frac{15}{\rho_i^4}A_i -\frac3{\rho_i^2}(1-\frac{\delta}2 ) \right] - \\
  &  & \hspace{-5mm} - \; \frac{12}{\rho_i^2} \left[\frac{\pd h_j(x)}{\pd x_k}\|x-M_i\| +h_j(x)
  \frac{(x_k-m^i_k)}{\|x-M_i\|} \right] - \nonumber \\
  &  & \hspace{-5mm} - \; \frac{36}{\rho_i^2}h_k(x) \frac{(x_j-m^i_j)}{\|x-M_i\|} + 3 (x_j-m^i_j)\frac{(x_k-m^i_k)}{\|x-M_i\|} \times \nonumber \\
  &  & \hspace{-5mm} \times \left[ -\frac{12}{\rho_i^2} \frac{<\!x-M_i,T-M_i\!>}{\| x-M_i \|} +
  \frac{10}{\rho_i^3}A_i +\frac3{\rho_i} (1-\frac{\delta}2 ) \right], \nonumber
\end{eqnarray}
with
$$
 \frac{\pd h_j(x)}{\pd x_k}=(T_j-m^i_j)\frac{(x_k-m^i_k)}{\|x-M_i\|} -
 \frac{h_k(x)}{\|x-M_i\|^2}(x_j-m^i_j),
$$
and
$$
 h_k(x)=(T_k-m^i_k) \|x-M_i\| -<\!x-M_i,T-M_i\!>
 \frac{(x_k-m^i_k)}{\|x-M_i\|},
$$
\vspace{6pt} \noindent while for pure partial derivatives $\pd ^2
f(x)/\pd x_j^2$ we have
 \be \label{d2(jj)f_quintic}
  \frac{\pd ^2f(x)}{\pd x_j^2}= \left\{
  \begin{array}{ll}
    \frac{\pd ^2Q_i (x)}{\pd x_j^2}, & x\in S_i, \, i \in \{2,\ldots,m\} , \\
     2, & x \notin S_2 \cup \ldots \cup S_m
     \, ,
   \end{array} \right.
 \ee
for $j=1,\ldots, N,$ and
 \begin{eqnarray} \label{d2(jj)Q_quintic}
  \frac{\pd^2 Q_i (x)}{\pd x_j^2} & = & -\frac6{\rho_i^4}
  \left[\frac{\pd h_j(x)}{\pd x_j} \|x-M_i\|^3+3h_j(x)(x_j-m^i_j)
  \|x-M_i\| \right] + \nonumber \\
  &  & \hspace{-5mm} + \; \left[ 5\|x-M_i\|^3+15(x_j-m^i_j)^2\|x-M_i\| \right] \times \nonumber \\
  &  & \hspace{-5mm} \times \left[ -\frac6{\rho_i^4} \frac{<\!x-M_i,T-M_i\!>}{\| x-M_i \|}+
  \frac6{\rho_i^5}A_i +\frac1{\rho_i^3} (1-\frac{\delta}2 ) \right] - \nonumber \\
  &  & \hspace{-5mm} - \; \frac{30}{\rho_i^4}h_j(x)(x_j-m^i_j)
  \|x-M_i\| + \\
  &  & \hspace{-5mm} + \; \frac{16}{\rho_i^3} \left[\frac{\pd h_j(x)}{\pd x_j} \|x-M_i\|^2+2h_j(x)
  (x_j-m^i_j) \right]  + \nonumber \\
  &  & \hspace{-5mm} + \; \frac{64}{\rho_i^3}h_j(x)(x_j-m^i_j) + \left[
  4\|x-M_i\|^2 +8(x_j-m^i_j)^2 \right]\times \nonumber \\
  &  & \hspace{-5mm} \times \left[ \frac{16}{\rho_i^3} \frac{<\!x-M_i,T-M_i\!>}{\| x-M_i\|}-
  \frac{15}{\rho_i^4}A_i -\frac3{\rho_i^2} (1-\frac{\delta}2 ) \right] -  \nonumber \\
  &  & \hspace{-5mm} - \; \frac{12}{\rho_i^2} \left[\frac{\pd h_j(x)}
  {\pd x_j} \|x-M_i\| +h_j(x) \frac{(x_j-m^i_j)}{\|x-M_i\|} \right] - \nonumber\\
  &  & \hspace{-5mm} - \; \frac{36}{\rho_i^2}h_j(x) \frac{(x_j-m^i_j)}{\|x-M_i\|} + \left[
  3\|x-M_i\|+3\frac{(x_j-m^i_j)^2}{\|x-M_i\|} \right]\times \nonumber \\
  &  & \hspace{-5mm} \times \left[ -\frac{12}{\rho_i^2} \frac{<\!x-M_i,T-M_i\!>}{\| x-M_i\|}+
  \frac{10}{\rho_i^3}A_i +\frac3{\rho_i} (1-\frac{\delta}2 ) \right] + \delta, \nonumber
\end{eqnarray}
with
$$
 \frac{\pd h_j(x)}{\pd x_j}=(T_j-m^i_j)\frac{(x_j-m^i_j)}{\|x-M_i\|} -
 \frac{h_j(x)}{\|x-M_i\|^2}(x_j-m^i_j)-\frac{<\!x-M_i,T-M_i\!>}{\|x-M_i\|}.
$$

\medskip

\bibliographystyle{plain}
\bibliography{Test_Gen_arxiv}



\end{document}